# Buffon's Triangle - A Variant of the Buffon Needle Method for a Probabilistic Determination of the Value of π

Devlin Gualtieri
Tikalon LLC, Ledgewood, New Jersey
(gualtieri@ieee.org)

*I present a variant of the Buffon Needle method for determination of the value of the mathematical constant, π. The original method is based on the random casting of a needle of length ℓ onto a planked floor of plank width **L**. The described variant involves the random casting of an equilateral triangle with side length ℓ onto a tiled floor consisting of square tiles of side width **L**. Source code for the computer simulation of this method is provided.*

Introduction

I present a variant of the Buffon Needle method for determination of the value of the mathematical constant, π, devised in 1777 by the French mathematician, Georges-Louis Leclerc, Comte de Buffon.[1,2] The Buffon needle method is based on the random casting of a needle of length ℓ onto a planked floor of plank width **L**. The described variant involves the random casting of an equilateral triangle with side length ℓ onto a tiled floor consisting of square tiles of side width **L**.

Buffon's estimate for π in the needle problem for the case where ℓ is not greater than **L** is given by the simple equation,

$$\pi = (2/p)(\ell/L) \qquad (1)$$

for which the probability **p** is is the probability that the needle will lie across the seam between adjacent floor planks. The usual case is when the needle length equals the distance between floor planks; i. e, (ℓ/L) = 1. I used this case in my Buffon's triangle model.

Buffon's Triangle Model

The model appears in schematic form in fig. 1.

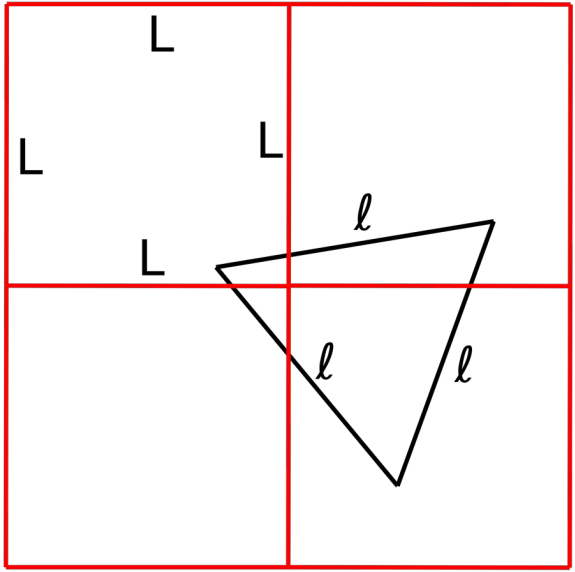

Fig. 1. Schematic diagram of the Buffon's Triangle model. The side lengths of the equilateral triangle are equal to the side lengths of the tile squares.

The estimate for π in this case where the side length of the equilateral triangle $\ell$ is equal to the side length of each square of the tiling **L**, is given by

$$\pi = 12(\text{trials}/\text{intersections}) \qquad (2)$$

for which **trials** is the number of times that the triangle is randomly cast onto the square tiling, and **intersections** is a count of the number of times that a side of the triangle intersects a seam of the square tiling.

Computer simulation

The source code for the computer simulation of Buffon's Triangle appears as Appendix I. This program includes an option for creating images of a few random triangle casts onto the square tiling. Fig. 2 contains examples of these castings.

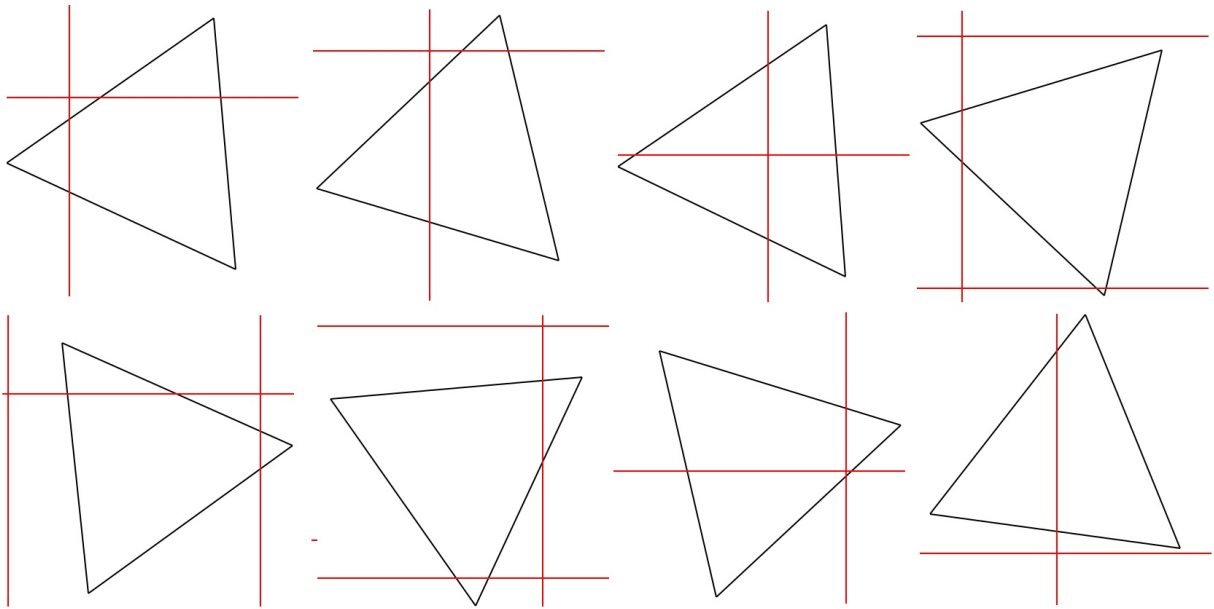

Fig. 2. Examples of equilateral triangle castings onto the square tile array for the Buffon's Triangle model.

A thousand runs of a million trials generates a histogram with a mean value of 3.14157, as shown in fig. 3.

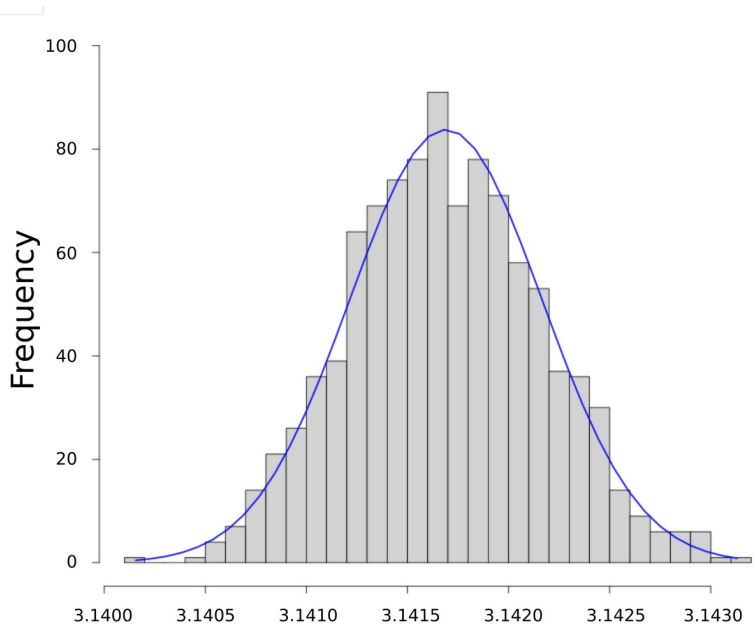

Fig. 3. Histogram of the estimated value of pi in the Buffon's Triangle model. These data are for 1000 runs of a million trials each.

# References

1. Histoire de l'Acad. Roy. des. Sciences, vol. 43-45; Histoire naturelle, générale et particulière Supplément 4 (1777), pp. 46ff. (via Google Books)
2. Buffon's needle problem, Wikipedia

# Appendix I – Source code for buffon.c.

```c
/* -*- Mode: C; indent-tabs-mode: t; c-basic-offset: 4; tab-width: 4 -*- */
/*
 * buffon.c
 * Copyright (C) DMGualtieri 2024 <gualtieri@ieee.org>
 *
 * buffon is free software: you can redistribute it and/or modify it
 * under the terms of the GNU General Public License as published by the
 * Free Software Foundation, either version 3 of the License, or
 * (at your option) any later version.
 *
 * buffon is distributed in the hope that it will be useful, but
 * WITHOUT ANY WARRANTY; without even the implied warranty of
 * MERCHANTABILITY or FITNESS FOR A PARTICULAR PURPOSE.
 * See the GNU General Public License for more details.
 *
 * You should have received a copy of the GNU General Public License along
 * with this program.  If not, see <http://www.gnu.org/licenses/>.
 */

// compile using gcc: gcc -o buffon buffon.c -lm

#include <stdio.h>
#include <stdlib.h>
#include <time.h>
#include <math.h>
#include <string.h>

#define trials 1000000
#define number_of_images 20
#define create_images 1            //1 to create example images, 0 to not

/* Prototypes */
char *strcpy(char *dest, const char *src);
void srand(unsigned seed);
int random_coordinate(int m);
float random_angle(void);
/* end of prototypes */

int i, j, k, n;
int r, x, y, x_0, y_0, x_1, y_1, x_2, y_2;
int x_line, y_line;        //lines that define x-y grid
int apex_x[3];             //apex x coordinates sorted from low to high
int apex_y[3];             //apex y coordinates sorted from low to high
int temp;                  //temporary storage for sorting
int count_x = 0;
int count_y = 0;
float angle;
char str_temp[8] = "";
char svg_filename[32] = "plot_000.svg";
FILE *svg_file;

//generates a psuedo-random integer between 0 and m
int random_coordinate(int m)
{
    return rand() % m;
}

//generates a psuedo-random float between 0.0 and 2pi radians
float random_angle(void)
{
    return (float) 6.28319 *((double) rand()) / RAND_MAX;
}
```

```c
int main(void)
{

/* seed random number function */
    srand((unsigned) time((time_t *) NULL));

//using an equilateral triangle inscribed in a circle
//which allows easy rotation.  120 degrees = 2.0944 radians
//x = r cos(angle), y = r sin(angle).  radius set at 10000
//and we add 10000 to keep everything positive
    r = 10000;

    for (i = 0; i < trials; i++) {
        angle = random_angle();
        x_0 = 10000 + (r * cos(angle));
        apex_x[0] = x_0;
        y_0 = 10000 + (r * sin(angle));
        apex_y[0] = y_0;
        angle = angle + 2.0944;
        x_1 = 10000 + (r * cos(angle));
        apex_x[1] = x_1;
        y_1 = 10000 + (r * sin(angle));
        apex_y[1] = y_1;
        angle = angle + 2.0944;
        x_2 = 10000 + (r * cos(angle));
        apex_x[2] = x_2;
        y_2 = 10000 + (r * sin(angle));
        apex_y[2] = y_2;

//sort apex x and y coordinates

        for (j = 0; j < 3; j++) {
            for (k = 0; k < 2; k++) {
                if (apex_x[k] > apex_x[k + 1]) {
                    temp = apex_x[k + 1];
                    apex_x[k + 1] = apex_x[k];
                    apex_x[k] = temp;
                }
            }
        }

        for (j = 0; j < 3; j++) {
            for (k = 0; k < 2; k++) {
                if (apex_y[k] > apex_y[k + 1]) {
                    temp = apex_y[k + 1];
                    apex_y[k + 1] = apex_y[k];
                    apex_y[k] = temp;
                }
            }
        }

//triangle sides are 17321
        x_line = random_coordinate(17321);
        y_line = random_coordinate(17321);

//test for intersections of x grid lines with triangle segments

        for (n = x_line; n < 20000; n = n + 17321) {
            if ((n > apex_x[0]) && (n <= apex_x[1])) {
                count_x++;
            }

            if ((n > apex_x[0]) && (n <= apex_x[2])) {
                count_x++;
            }

            if ((n > apex_x[1]) && (n <= apex_x[2])) {
                count_x++;
            }
        }

        for (n = y_line; n < 20000; n = n + 17321) {
            if ((n > apex_y[0]) && (n <= apex_y[1])) {
                count_y++;
            }

            if ((n > apex_y[0]) && (n <= apex_y[2])) {
                count_y++;
            }
```

```c
            if ((n > apex_y[1]) && (n <= apex_y[2])) {
                count_y++;
            }
        }

        if ((create_images == 1) && ((i < number_of_images))) {
//create SVG file name
            sprintf(str_temp, "%02d", i);
            strcpy(svg_filename, "plot");
            strcat(svg_filename, str_temp);
            strcat(svg_filename, ".svg");

            if ((svg_file = fopen(svg_filename, "w")) == NULL) {
                printf("\nSVG output file cannot be opened.\n");
                exit(1);
            }

            // Print svg header
            fprintf(svg_file,
                    "<svg height=\"400\" width=\"400\" xmlns=\"http://www.w3.org/2000/svg\">\n");
            fprintf(svg_file,
                    "<rect x=\"0\" y=\"0\" width=\"400\" height=\"400\" fill=\"white\" />\n");

//print svg data
            fprintf(svg_file,
                    "<line x1=\"%d\" y1=\"%d\" x2=\"%d\" y2=\"%d\" ",
                    x_0 / 50, y_0 / 50, x_1 / 50, y_1 / 50);
            fprintf(svg_file, "stroke=\"black\" stroke-width=\"2\"/>\n");
            fprintf(svg_file,
                    "<line x1=\"%d\" y1=\"%d\" x2=\"%d\" y2=\"%d\" ",
                    x_1 / 50, y_1 / 50, x_2 / 50, y_2 / 50);
            fprintf(svg_file, "stroke=\"black\" stroke-width=\"2\"/>\n");
            fprintf(svg_file,
                    "<line x1=\"%d\" y1=\"%d\" x2=\"%d\" y2=\"%d\" ",
                    x_2 / 50, y_2 / 50, x_0 / 50, y_0 / 50);
            fprintf(svg_file, "stroke=\"black\" stroke-width=\"2\"/>\n");
            fprintf(svg_file,
                    "<line x1=\"%d\" y1=\"%d\" x2=\"%d\" y2=\"%d\" ",
                    x_line / 50, 0, x_line / 50, 400);
            fprintf(svg_file, "stroke=\"red\" stroke-width=\"2\"/>\n");
            fprintf(svg_file,
                    "<line x1=\"%d\" y1=\"%d\" x2=\"%d\" y2=\"%d\" ",
                    0, y_line / 50, 400, y_line / 50);
            fprintf(svg_file, "stroke=\"red\" stroke-width=\"2\"/>\n");
            fprintf(svg_file,
                    "<line x1=\"%d\" y1=\"%d\" x2=\"%d\" y2=\"%d\" ",
                    (x_line + 17321) / 50, 0, (x_line + 17321) / 50, 400);
            fprintf(svg_file, "stroke=\"red\" stroke-width=\"2\"/>\n");
            fprintf(svg_file,
                    "<line x1=\"%d\" y1=\"%d\" x2=\"%d\" y2=\"%d\" ",
                    0, (y_line + 17321) / 50, 400, (y_line + 17321) / 50);
            fprintf(svg_file, "stroke=\"red\" stroke-width=\"2\"/>\n");

// Print svg footer
            fprintf(svg_file, "</svg>");

            fclose(svg_file);

        }

    }

    printf("\ncount_x = %d\tcount_x/trials = %f", count_x,
           (float) count_x / trials);
    printf("\ncount_y = %d\tcount_y/trials = %f\n", count_y,
           (float) count_y / trials);
    printf("\npi estimate = %f\n",
           ((float) 12 * trials / (count_x + count_y)));

    printf("Done.\n");

    return 0;

}
```